\documentclass[epic,eepic,11pt]{amsart}

\usepackage{amsfonts,mathrsfs}
\usepackage[mathscr]{eucal}

\usepackage{amssymb}
\usepackage{amscd}
\usepackage{fancyhdr}
\usepackage{euler,eucal}

\DeclareMathAlphabet{\mathbf}{T1}{ppl}{bx}{n}
\DeclareMathAlphabet{\mathrm}{T1}{ppl}{m}{n}

\numberwithin{equation}{section}

\newcommand\note[1]%
{$^\dagger$\marginpar{\footnotesize{$^\dagger${#1}}}}

\def\({\left(}
\def\){\right)}
\def\<{\left<}
\def\>{\right>}

\def\newi {\sqrt{-1}\, }

\newtheorem{theorem}{Theorem}[section]
\newtheorem{proposition}[theorem]{Proposition}
\newtheorem{lemma}[theorem]{Lemma}
\newtheorem{definition}[theorem]{Definition}

\theoremstyle{definition}
\newtheorem{example}[theorem]{Example}
\newtheorem{remark}[theorem]{Remark}

\newcommand\lie{\mathfrak}
\renewcommand\t{\lie{t}}
\newcommand\g{\lie{g}}
\newcommand\fk{\lie{k}}

\newcommand\R{\mathbb{R}}
\newcommand\C{\mathbb{C}}

\renewcommand\P{\mathbb{P}}
\newcommand\CP{\mathbb{CP}}

\newcommand\J{\mathcal{J}}
\newcommand\cO{\mathcal{O}}

\newcommand     {\comment}[1]   {}
\newcommand{\mute}[2] {}
\newcommand     {\printname}[1] {}

\newcommand\func[1]{\operatorname{\mathrm{#1}}}
\newcommand\funclim[1]{\operatorname*{\mathrm{#1}}}

\newcommand\type{\func{type}}

\renewcommand\dim{\func{dim}}

\renewcommand\lim{\funclim{lim}}

\newcommand\sur{\mathrel{\to\kern-1.8ex\to}}
\newcommand\iso{\mathrel{\hookrightarrow\kern-1.8ex\to}}

\newcommand\longhookrightarrow{\lhook\joinrel\longrightarrow}

\newcommand\longsur{\mathrel{\longrightarrow\kern-1.8ex\to}}
\newcommand\longiso{\mathrel{\longhookrightarrow\kern-1.8ex\to}}

\begin{document}

\bibliographystyle{amsalpha}
\date{\today}

\title{Symmetries in generalized K\"{a}hler geometry}

\author{Yi Lin, \,\,\,  Susan Tolman}
\address{Department of Mathematics, University
of Toronto, Canada, M5S2E4, yilin@math.toronto.edu}
\address{Department of Mathematics, University of Illinois at
Urbana-Champaign, Urbana, IL, USA, 61801, stolman@math.uiuc.edu}

\begin{abstract}
We define the notion of a moment map and reduction in both
generalized complex geometry and generalized K\"{a}hler geometry. As
an application, we give very simple explicit constructions of
bi-Hermitian structures on $\C\P^n$, Hirzebruch surfaces, the blow
up of $\CP^2$ at arbitrarily many points, and other toric varieties,
as well as complex Grassmannians.

\end{abstract}

\maketitle

\section{Introduction}

Generalized complex structures were introduced by N. Hitchin
\cite{H02}, and further developed by Gualtierie \cite{Gua}. It
contains both symplectic and complex structures as extremal special
cases, and provides a useful differential geometric context for
understanding  some recent development in string theory. An
associated notion of generalized K\"{a}hler structure was introduced
by Gualtieri \cite{Gua}, who shows that this notion is essentially
equivalent to that of a bi-Hermitian structure, which was first
discovered by physicists (see \cite{GHR84} ) studying
super-symmetric nonlinear $\sigma$-models.

The theory of bi-Hermitian geometry suffered from a lack of
interesting examples. As stated in \cite{AGG99} (see also
\cite{AGG05}), an important open problem in this field was to
determine whether or not there exist bi-Hermitian structures on
$\C\P^2$, a minimal ruled surface admitting an effective
anti-canonical divisor, or a complex surface obtained from them by
blowing up points along an effective anti-canonical divisor.

Using the deformation theorem he developed for generalized complex
structure, Gualtieri proved \cite{Gua} that there exists a
bi-Hermitian structure on $\CP^2$. More recently, Hitchin \cite{H05}
used the generalized K\"{a}hler geometric approach to give an
explicit construction of a bi-Hermitian structure on $\CP^2$ and
also on $\CP^1 \times \CP^1$.

For manifolds with symmetries, the related notions of moment maps
and quotient are important in many  geometries. It is an interesting
question if there exist natural notions of moment maps and quotients
in generalized complex and K\"{a}hler geometries. Some attempts have
been made in this direction.  In \cite{Cr04}, Cranic proposed a
definition of moment maps in generalized complex geometry. But it
appears that the condition to make his definition work is rather
restrictive. In \cite{H05} Hitchin also presented a quotient
construction in generalized K\"{a}hler geometry which works for
certain interesting special cases. More recently, we have discovered
that several other groups are independently working on related
projects, including: Bursztyn, Cavalcanti, and Gualtieri; and Xu
Ping.

In this paper, we define the notions of generalized moment map for a
compact Lie group action on a generalized complex  manifold. Using
this definition we then define a generalized complex structure on
the reduced space, which is natural up to transformation by an exact
$B$-field, that is, the reduced space  has a natural  equivalence
class of generalized complex structures in the sense specified in
\cite{H02}. Moreover, we show that the quotient structure has the
same type as the original generalized complex structure. In the case
that the generalized complex structure is derived from a symplectic
structure, our definitions agree with the standard definitions of
moment map and symplectic reduction. Compared with the definition of
moment maps given in \cite{Cr04}, our approach works in greater
generality.

We then consider the compact Lie group action on a generalized
K\"{a}hler manifolds; in this case, the generalized moment map is
simply the generalized moment map for the first generalized complex
structure. Finally, we define a natural generalized  K\"{a}hler
structure on the reduced space, and give formulas for the types of
this structure. Again, in the case that the generalized K\"{a}hler
structure is derived from a K\"{a}hler structure, this agrees with
the usual K\"{a}hler reduction.

As an application, we give a very simple explicit construction of
bi-Hermitian structures on
 $\C\P^n$, Hirzebruch surfaces, the blow up of $\CP^2$ at arbitrarily many
points, and other toric varieties, as well as complex Grassmannians.
As shown in this paper, in practice our method will give one a
powerful machinery of producing bi-Hermitian structures on any
manifold which can be produced as the symplectic quotient of $\C^N$.

The plan of this paper is as follows.

Section 2 shows that under reasonable assumptions the Courant
bracket is preserved under  restriction and quotient.

Section 3 defines generalized moment map for a compact Lie group
acting on a generalized complex manifold, and constructs a
generalized complex structure on the reduced space at every regular
value.

Section 4  extends the results of Section 3 to the generalized
K\"{a}hler case. We also discuss the connection between the
hyper-K\"{a}hler quotient and the generalized K\"{a}hler quotient.

Section 5 presents the explicit constructions of bi-Hermitian
structures.

\section{The Courant bracket}

Let $V$ be a $n$ dimensional vector space. There is a natural metric
of type $(n,n)$ on $V \oplus V^*$ given by
$$\langle X + \alpha, Y + \beta \rangle = \frac{1}{2}
( \alpha(Y) + \beta(Y)).$$ Given a subspace $F \subset V \oplus
V^*$, let $F^\perp \subset V \oplus V^*$ denote the perpendicular
with respect to this metric. In contrast, if $F \subset V$ (or
$V^*$), let $F^0 \subset V^*$ (or $V$) denote the annihilator of
$F$.
%%S added
Moreover, let $\pi: V_\C \oplus V_\C^* \to V_\C$ denote the natural
projection.

Let $M$ be a $n$ dimensional manifold. There is a natural metric of
type $(n,n)$ on $TM \oplus T^*M$ given by
$$\langle X + \alpha, Y + \beta \rangle = \frac{1}{2}
( \alpha(Y) + \beta(Y)),$$ which extends naturally to $T_\C M\oplus
T_\C ^*M=(TM\oplus T^*M)\otimes \C$.

%%S added both cases.
Given a subbundle $F \subset T M\oplus T ^*M$ (or $T_\C M \oplus
T_\C^* M)$, let  $F^\perp \subset T M\oplus T ^*M$ or $(T_\C M
\oplus T_\C^* M)$ denote the perpendicular of $F$ with respect to
the above metric.
%%S added
Moreover, let $\pi: T_\C M \oplus T^*_\C M  \to T_\C M$ denote the
natural projection.

The {\bf Courant bracket} on $T_\C M \oplus T_\C ^*M$ is defined by
$$ [X + \alpha, Y + \beta]  = [X,Y] + L_X \beta - L_Y \alpha
- \frac{1}{2} (d \iota_X \beta - d \iota_Y \alpha). $$

We will now examine how the Courant bracket behaves under
restriction and quotient.

%%S changed $df$ to be real.
\begin{lemma} \label{le:closed1}
Let $M$ be a manifold and let $\g^*$ be the dual of a vector space
$\g$. Given a submersion  $f \colon M \to \g^*$, let $df \subset T M
\oplus T ^*M$ denote the subbundle spanned by the differentials
$df^{\xi}$ for $\xi \in \g$. Then $df_\C ^\perp$ is closed under the
Courant bracket.

Moreover, the restriction from $M$ to $f^{-1}(0)$ induces a natural
map from $df_\C^\perp \subset T_\C M \oplus T_\C^*M$ to
$T_\C(f^{-1}(0) ) \oplus T_\C^*(f^{-1}(0))$. If $\Gamma$ is a
sub-bundle of $T_\C M\oplus T_\C ^*M$ which is closed under the
Courant bracket, then the image of $\Gamma \cap df_\C^\perp$ under
this map is also closed under the Courant bracket.
\end{lemma}

\begin{proof} Let $X+\alpha$ and  $Y+\beta$ be sections of
$df_\C^{\perp}$. Given any $\xi \in g$, by assumption $\iota_X df
\xi=\iota_Y df^\xi=0$. Hence by Cartan's Formula
$$ \iota_{[X,Y]} df^\xi = L_X \iota_Y df^\xi - L_Y \iota_X df^\xi +
\iota_X \iota_Y d df^\xi - d\iota_X \iota_Y df^\xi = 0.$$ Therefore,
$[X + \alpha, Y + \beta] \in df_\C^{\perp}$. This proves the first
claim.

Finally, if $\Gamma$ is closed under the Courant bracket, then
$\Gamma \cap df_\C^\perp$ is also closed. Since a straightforward
check of the definition shows that map induced by restriction
preserves the Courant bracket, the second claim is obvious.
\end{proof}

\begin{lemma}\label{le:closed2} Let a compact Lie group
$G$ act freely on a manifold $M$, and let $\g_M \subset T_\C M
\oplus T_\C ^*M $ denote the subbundle spanned by the fundamental
vector fields $\xi_M$ for $\xi$ is in the Lie algebra $\g$ of $G$.
Then the set of $G$-invariant sections of $(\g_M)_\C^\perp$ is also
closed under the Courant bracket.

Moreover,  the quotient map from $M$ to $M/G$ induces a natural map
from the set of $G$-invariant sections of $(\g_M)_\C^\perp \subset
T_\C M \oplus T_\C ^*M$ to the section of $T_\C (M/G) \oplus T_\C
^*(M/G)$. Let $\Gamma$ be an $G$-invariant sub-bundle of $T_\C M
\oplus T_\C ^*M$ which is closed under the Courant bracket. Then the
image of $\Gamma \cap (g_M)_\C^{\perp}$ under this map is also
closed under the Courant bracket.
\end{lemma}

\begin{proof}
Let $X+\alpha$ and  $Y+\beta$ be $G$-invariant sections of
$(g_M)_\C^{\perp}$. Given any $\xi \in \g$, by assumption
$\iota_{\xi_M} \alpha=\iota_{\xi_M} \beta =0$. Since $X+\alpha$ and
$Y + \beta$ are $G$ invariant, $[\xi_M,X] = 0$ and $L_{\xi_M}
\iota_X \beta = 0$. Therefore
$$\iota_{\xi_M} L_X
\beta = \iota_{[\xi_{M},X]} \beta  + L_X \iota_{\xi_M} \beta = 0
\quad \mbox{and} \quad \iota_{\xi_M} d \iota_X \beta = L_{\xi_M}
\iota_X \beta + d \iota_X\iota_{\xi_M} \beta = 0.$$ Similarly,
$\iota_{\xi_M} L_Y \alpha = \iota_{\xi_M} d \iota_Y \alpha = 0$.
Hence, $[X+\alpha, Y+\beta] \in (g_M)_\C^{\perp}$.

Finally, if $\Gamma$ is closed under the Courant bracket, then
$\Gamma \cap (\g_M)_\C^\perp$ is also closed. Since a
straightforward check of the defintion shows that the map induced by
the quotient map preserves the Courant bracket, the second claim is
obvious.
\end{proof}

\section{Generalized complex structures}
\label{hamiltonian gcmfld}

A {\bf generalized complex structure} on a vector space $V$ is an
orthogonal linear map $\mathcal{J} \colon V \oplus V^* \to V \oplus
V^*$ so that $\mathcal{J}^2 = -1$. Given a generalized complex
structure $\mathcal{J}$, let $L \subset  V_\C \oplus V_\C^*$ be the
$\sqrt{-1}$ eigenspace of $\mathcal{J}$. Then $L$ is maximal
isotropic and $L \cap \overline{L} = \{0\}$. Conversely, given a
maximal isotropic $L \subset V_\C \oplus V_\C^*$ so that $L \cap
\overline{L} = \{0\}$, there exists a unique generalized complex
structure on $V$ whose $\sqrt{-1}$ eigenspace is $L$.

Let $\mathcal{J}$ be a generalized complex structure on a vector
space $V$ and let $W = V \oplus V^*$. If $P \subset W$ is a
$\mathcal{J}$-invariant subspace, then since $\mathcal{J}$ is
orthogonal there is a  restriction map $\overline{\mathcal{J}}
\colon P^{\perp} \to P^{\perp}$. If $P$ is also isotropic, let
$\widetilde{W}= P^{\perp}/P$; there is a  quotient map
$\widetilde{\mathcal{J}} \colon \widetilde{W} \to \widetilde{W}$.
Clearly, $\overline{\mathcal{J}}^2 = -1$,
$\widetilde{\mathcal{J}}^2 = -1$ and both maps are orthogonal.
%%S NEW MATH
Also, if $L$ is the $\sqrt{-1}$ eigenbundle of $J$, then $L \cap
P_\C^\perp$ is the $\sqrt{-1}$ eigenbundle of $\overline{\J}$, and
the image of $L \cap P_\C^\perp$ in $\widetilde{W}_\C$ is  the
$\sqrt{-1}$ eigenbundle of $\widetilde{\J}$. Finally, if $P = (P
\cap V) \oplus (P \cap V^*) $, let $\widetilde{V}$ be the quotient
of $(P \cap V^*)^0 \subset V$ by $(P \cap V)$; the spaces
$\widetilde{W}$ and $\widetilde{V} \oplus \widetilde{V}^*$ are
naturally isomorphic. Hence $\mathcal{J}$ naturally induces a
generalized complex structure $\widetilde{\mathcal{J}}$ on
$\widetilde{V}$.

The {\bf type} of $\mathcal{J}$ is the codimension of $\pi(L)$ in
$V_\C$, where $L$ is the $\sqrt{-1}$ eigenspace of $\mathcal{J}$.
Recall that $\pi \colon V_\C \oplus V^*_\C  \to V_\C$ is the natural
projection.) The following lemma will help us compute types.

\begin{lemma}\label{pretype}
Let $\mathcal{J}$ be a  generalized complex structure on a vector
space $V$, and let $L \subset V_\C \oplus V^*_\C$ be its $\sqrt{-1}$
eigenspace.
 If a subspace $R \subset V_\C \oplus V^*_\C$ satisfies
$\mathcal{J}(R) \cap R = \{0\}$, then
$$  \dim (\pi(L \cap R^\perp \cap \mathcal{J}(R)^\perp))
=  \dim(\pi(L + R)) - \dim(R).$$
\end{lemma}

\begin{proof}
Since $L$ is the $\sqrt{-1}$ eigenspace of $\mathcal{J}$,
$$L \cap R^\perp \cap \mathcal{J}(R)^\perp = L \cap R^\perp.$$
Since $L \cap R^\perp \cap \mathcal{J}(R)^\perp$ is the $\sqrt{-1}$
eigenspace of the restriction of $\mathcal{J}$ to $R^\perp \cap
\mathcal{J}(R)^\perp$ and $\J(R) \cap R = \{0\}$,
$$\dim(L \cap R^\perp \cap \mathcal{J}(R)^\perp) = \dim V - \dim R.$$
Since $V_\C^*$ is the kernel of $\pi$,
$$\dim( \pi(L \cap R^\perp)) =
\dim (L \cap R^\perp) - \dim(L \cap R^\perp \cap V_\C^*).$$ Finally,
since $L$ is maximal isotropic, $L = L^\perp$, and so $$L \cap
R^\perp \cap V_\C^* = (L + R)^\perp \cap V_\C^* = \pi(L + R)^0.$$
\end{proof}

\begin{lemma} \label{type}
Let $\mathcal{J}$ be a generalized complex structure on a vector
space $V$. Consider a subspace $Q \subset V$ so that $\mathcal{J}(Q)
\subset V^*$ and so that $P= Q \oplus \mathcal{J}(Q) \subset V
\oplus V^*$ is isotropic. Let $\widetilde{\mathcal{J}}$  be the
natural generalized complex structure on $\widetilde{V} =
\mathcal{J}(Q)^0/Q$. Then
$$\type(\widetilde{\mathcal{J}})= \type(\mathcal{J}).$$
\end{lemma}

\begin{proof}
Let $L$ and $\widetilde{L}$ be the $\sqrt{-1}$ eigenspaces of
$\mathcal{J}$ and $\widetilde{\mathcal{J}}$, respectively; let $\pi
\colon V_\C \oplus V^*_\C \to V_\C$ and $\widetilde{\pi} \colon
\widetilde{V}_\C \oplus \widetilde{V}_\C^* \to \widetilde{V}_\C$ be
the natural projections.

Since $Q \subset V$ and $\mathcal{J}(Q) \subset V^*$, it is
immediately clear that $Q_\C \subset \pi(L)$ and $Q \cap
\mathcal{J}(Q) = \{0\}$. Therefore, by Lemma~\ref{pretype} $\dim
(\pi(L \cap Q_\C^\perp \cap \mathcal{J}(Q_\C)^\perp)) = \dim(\pi(L))
- \dim(Q). $ Moreover, $\widetilde{\pi}(\widetilde{L}) $ is the
projection of $\pi(L \cap Q_\C^{\perp} \cap \mathcal{J}(Q_\C)^\perp
) \subset  \mathcal{J}(Q)_\C^0$ to $\widetilde{V}_\C =
\mathcal{J}(Q)_\C^0/Q_\C$, so $\dim(\widetilde{\pi}(\widetilde{L}))
= \dim(\pi(L \cap Q_\C^{\perp} \cap \mathcal{J}(Q_\C)^\perp)) -
\dim(Q)$. Finally, $\dim(\widetilde{V}) = \dim(V) - 2 \dim (Q)$.
\end{proof}

A {\bf generalized almost complex structure} on a manifold $M$ is an
orthogonal  bundle map $\mathcal{J} \colon TM \oplus T^*M \to TM
\oplus T^*M$ so that $\mathcal{J}^2 = -1$. Moreover, $\mathcal{J}$
is a {\bf generalized complex structure} if the $\sqrt{-1}$
eigenbundle of $\mathcal{J}$, $L \subset T_\C M \oplus T_\C^* M$, is
closed under the Courant bracket. The {\bf type} of $\mathcal{J}$ at
$m \in M$ is the type of the restricted generalized complex
structure on $T_m M$.

We now introduce several standard examples, as described in
\cite{Gua}.

\begin{example}\label{ex1}(\cite{Gua})
\begin{itemize}
\item[a)]
Let $(M, \omega)$ be a symplectic manifold. Then
$$ \mathcal{J}_{\omega}=\left( \begin{matrix} 0 & \omega^{-1}
\\ -\omega & 0 \end{matrix} \right)$$ is a generalized complex
structure of type $0$; the $\sqrt{-1}$ eigenspace of $\J_\omega$ is
$L_\omega = \{ X - \newi \iota_X \omega \mid X \in V_\C \}$.
\item[b)]
Let  $(M, J)$ be a $2n$ dimensional complex manifold. Then
$$ \mathcal{J}_J =\left( \begin{matrix}- J & 0 \\
0 & J^*  \end{matrix} \right)\,\,\,   $$ is a generalized complex
structure of type $n$; the $\sqrt{-1}$ eigenspace of $\J_J$ is $L_J
=  T_{0,1} \oplus T_{1,0}^*$, where $T_{1,0}$ is the $\sqrt{-1}$
eigenbundle of $J$.
\item[c)]
Let $(M_1, \mathcal{J}_1)$ and $(M_2, \mathcal{J}_2)$ be generalized
complex manifolds. Then
$$\mathcal{J}_1 \times \mathcal{J}_2 =\left(
\begin{matrix} \mathcal{J}_1 & 0
\\ 0 & \mathcal{J}_2 \end{matrix} \right).$$
is a  generalized complex structure on $M_1 \times M_2$, and
$$
\type(\mathcal{J}_1 \times \mathcal{J}_2)_{(m_1, m_2)} =
\type(\mathcal{J}_1)_{m_1} + \type(\mathcal{J}_2)_{m_2}.
$$
\item[d)] Let $B$ be a closed two-form on a manifold $M$, and
consider the  orthogonal bundle map $TM \oplus T^*M \to TM \oplus
T^*M$ defined by
$$e^B=\left( \begin{matrix} 1 & 0\\B & 1\end{matrix}\right),$$
where $B$ is regarded as a skew-symmetric map from $TM$ to $T^*M$.
This map preserves the Courant bracket (see \cite{Gua}.) As a simple
consequence, if $\mathcal{J}$ is a generalized complex structure on
$M$, then $\mathcal{J}'= e^B \mathcal{J} e^{-B}$ is another
generalized complex structure on $M$, called the {\bf
$\bold{B}$-transform} of $\mathcal{J}$. Moreover, the $\sqrt{-1}$
eigenbundle of $\J'$ is $e^B(L)$, so $\J$ and $\J'$ have the same
type. Finally, we say that $\mathcal{J}'$ is an {\bf  exact
$\bold{B}$-transform} if $B$ is exact.
\end{itemize}
\end{example}

%%S got rid of connect`
\begin{definition} \label{defmm}
Let a compact Lie group $G$ with Lie algebra $\g$ act on a manifold
$M$, preserving a generalized complex structure $\mathcal{J}$. Let
$L \subset T_\C M \oplus T_\C ^*M$ denote the $\sqrt{-1}$
eigenbundle of $\J$. A {\bf generalized moment map}  is a smooth
function $\mu  \colon  M \to \g_\C$ so that
\begin{itemize}
\item $\xi_M - \newi d\mu^\xi$ lies in $L$ for all $\xi \in \g$,
where $\xi_M$ denotes the induced vector field on $M$.
\item
$\mu$ is equivariant.
\end{itemize}
The generalized moment map $\mu$ is {\bf real} if  $\mu =
\overline{\mu}$. The action is {\bf Hamiltonian} if  a generalized
moment moment map exists.
\end{definition}

We are now ready to compute generalized moment maps for our basic
examples in Example~\ref{ex1}.

\begin{example}\label{ex2}\

\begin{itemize}

\item[a)]
Let $G$ acts on a symplectic manifold $(M, \omega)$ and let $\Phi
\colon M \to \g^*$ be a moment map, that is $\Phi$ is equivariant
and $ \iota_M \omega =  d\Phi^\xi$ for all $\xi \in \g$. Then $G$
also preserves $\mathcal{J}_\omega$, and  $\Phi$ is  also a (real)
generalized moment map for this action.

\item[b)] If $G$ acts on a
complex manifold $(M,J)$, preserving $J$, then $G$ also preserves
the associated generalized complex structure $\mathcal{J}_J$.
However, the action is never Hamiltonian because $\pi(L_J)$ contains
no real vectors.

\item[c)]
If $G$ acts on generalized complex manifolds $(M_1, \mathcal{J}_1)$
and $(M_2, \mathcal{J}_2)$ with generalized moment maps $\mu_1$ and
$\mu_2$, then the diagonal action of $G$ on the product manifold
$M_1 \times M_2$ preserves the generalized complex structure
$\mathcal{J}_1 \times \mathcal{J}_2$, and $\mu=\mu_1+\mu_2 \colon
M_1\times M_2 \to \g_\C$ is a generalized moment map for this
action.

\item[d)]
Let $G$ act on a generalized complex manifold $(M,\J)$ with
generalized moment map $\mu \colon  M \to \g_\C^*$. Let $B \in
\Omega^2(M)$ be closed and invariant, and let  $\Phi \colon M \to
\g^*$ be an equivariant map so that $\iota_{\xi_M} B = d \Phi^\xi$
for all $\xi \in \g$. Then $G$ also preserves the $B$-transform
$\J'$ of $\J$, and $\mu + \newi \Phi$ is a generalized moment map
for this action.
\end{itemize}
\end{example}

When the action is free, we can always perform an exact
$B$-transform so that the generalized moment map is real.

\begin{lemma}
\label{B-transform} Let a compact Lie group $G$  act freely on a
generalized complex manifold $(M, \mathcal{J})$ with generalized
moment map $\mu = f + \newi h \colon M \to \g_\C ^*$, where $f$ and
$h$ are real. Given a connection $\theta \in \Omega^1(M,\g)$, let $B
= d(\theta,h) \in \Omega^2(M;\R)$. The $B$-transform of $\J$ is an
invariant generalized complex structure $\mathcal{J}'$ with real
generalized moment map $f$. Here,  $(\theta,h) \in \Omega^1(M;\R)$
is obtained using the natural pairing on $\g$ and $\g^*$.

\end{lemma}

\begin{proof}
Since $\theta$ and $h$ are both equivariant, $(\theta,h)$  is
invariant. Since $B$ closed and invariant, $\mathcal{J}'$ an
invariant generalized complex structure. Since $\theta$ is a
connection and $(\theta,h)$ is invariant, $\iota_{\xi_M} B =
\iota_{\xi_M} d(\theta,h) = - d (\iota_{\xi_M} \theta,h) = - d
h^\xi$. The result now follows as in the example above.
\end{proof}

As shown in the next lemma, we get the second condition in the
definition of generalized moment map for free in the case of a torus
action.

\begin{lemma}\label{invariance}
Let a compact torus $T$  with Lie algebra $\t$ act on a manifold
$M$, preserving a generalized complex structure $\mathcal{J}$. Any
function $\mu \colon M \to \t_\C ^*$  which satisfies the first
criterion of Definition~\ref{defmm} is a generalized moment map.
\end{lemma}

\begin{proof} Fix $\xi$ and $\eta$  in $\t$.
Since $L$ is isotropic,
$$
0 =  \left< \xi_M-\newi d \mu^{\xi},\eta_M-\newi d \mu^{\eta}
\right> = -\newi \( \iota_{\xi_M} d \mu^\eta + \iota_{\eta_M} d
\mu^\xi\). $$ Since $T$ is abelian, $[\xi_M,\eta_M] = 0$. Since $d
\mu$ is closed, $L_{\xi_M} d \mu^\eta = d \iota_{\xi_M} d \mu^\eta$
and $L_{\eta_M} d \mu^\xi = d \iota_{\eta_M} d \mu^\xi$. Therefore,
\[\begin{split}
& [\xi_M-\newi d \mu^{\xi},\eta_M-\newi d \mu^{\eta}]\\
& =[\xi_M, \eta_M]- \newi L_{\xi_M}d\mu^{\eta} + \newi
L_{\eta_M}d\mu^{\xi} +
\dfrac{\newi}{2}(d\iota_{\xi_M}d\mu^{\eta}-d\iota_{\eta_M}d\mu^{\xi})\\
&= - \newi d\iota_{\xi_M}d\mu^{\eta}.
\end{split} \]
Since $L$ is closed under Courant bracket, this implies that
$d\iota_{\xi_M}d\mu^{\eta}\in L$.

Fix any orbit $\cO$ of the torus action. Since $L$ is isotropic, the
fact that $d\iota_{\xi_M}d\mu^{\eta}$ lies in $L$ implies that it
vanishes on $\cO$, and so $\iota_{\xi_M} d \mu^\eta$ is constant  on
$\cO$. Since $\cO$ is compact, this implies that
$\iota_{\xi_M}d\mu^{\eta}$ vanishes on $\cO$ Hence, $\mu^\eta$ is
invariant under the subgroup generated by $\xi$, and so $\mu$ is
$T$-invariant.
\end{proof}

\comment{
\begin{remark}
Gualtieri observed (see also \cite{AB04}) that for any generalized
complex manifold $(M,\J)$, the  bivector $\Pi$ defined by the upper
right quadrant of $\mathcal{J} \colon TM \oplus T^*M \rightarrow TM
\oplus T^*M$ is a real Poisson bivector; this gives rise to a
Poisson bracket $\{\cdot ,\cdot \}$ on $C^{\infty}(M)$. Let $\mu = f
+ \newi h$ be the generalized moment map of a Hamiltonian
generalized complex $G$-manifold, where $f$ and $h$ are real. A
simple calculation shows that the first condition in Definition
\ref{defmm} implies that $\Pi(df^{\xi})=-\xi_M$. Thus a necessary
condition for the existence of the generalized moment map is that
the action of the Lie group $G$ on the Poisson manifold $(M, \Pi)$
is Hamiltonian. And the second condition in the above definition
implies that the real moment map $f$ for the Poisson action is
equivariant, or equivalently, $\{f^\xi, f^\eta\}=f^{[\xi,\eta]}$ for
all $\xi$ and $\eta$ in $\g$. If $\mu' \colon M \to \g_\C$ is
another equivariant map, then an argument similar to the one given
in the proof of Lemma \ref{invariance} shows that $\mu'$ is a
generalized moment map if and only if $\mu'-\mu$ differs by a
constant on each symplectic leaf of the canonical Poisson structure
on $M$.
\end{remark}
}

Let a compact Lie group $G$ act on a generalized complex manifold
$(M,\mathcal{J})$ with generalized moment map $\mu = f + \newi h$,
where $f$ and $h$ are real. Let $\cO_a$ be the co-adjoint orbit
through $a \in \g^*$. Since $f$ is equivariant, $G$ acts on
$f^{-1}(\cO_a)$. Moreover, if $G$ acts freely on $f^{-1}(\cO_a)$,
then since $\mathcal{J}(df^\xi) = -\xi_M - dh^\xi$ for all $\xi \in
\g$, $\cO_a$ consists of regular values of $f$ and the {\bf
generalized complex quotient}
$$M_a =  f^{-1}(\cO_a)/G$$
is a manifold.

\begin{lemma} \label{complex reduction}
Let a compact Lie group $G$ act on a generalized complex manifold
$(M,\mathcal{J})$ with a real generalized moment map $f \colon M \to
\g_\C^*$. Let $\cO_a$ be the co-adjoint orbit through $a \in \g^*$.
If $G$ acts freely on $f^{-1}(\cO_a)$, the generalized complex
quotient $M_a$ inherits a natural generalized complex structure
$\widetilde{\J}$.

Moreover, for all $m \in f^{-1}(\cO_a)$,
$$\type(\widetilde{\mathcal{J}})_{[m]} = \type(\mathcal{J})_m.$$
\end{lemma}

\begin{proof}
First, assume that $a = 0$.

By restricting to a neighborhood of $f^{-1}(0)$, we may assume that
$G$ acts freely,  and that hence $f$ is a submersion. By the
definition of generalized moment map, $\mathcal{J}(\xi_M) = df^\xi$
for all $\xi \in \g$, so $\mathcal{J}(\g_M) = df$. Therefore, $\g_M
\oplus df$ is a $\mathcal{J}$-invariant subbundle of $TM \oplus
T^*M$. Since  $G$ acts on $f^{-1}(0)$, $\g_M \oplus df$ is also
isotropic when restricted to $f^{-1}(0)$. As in the discussion
preceding Lemma \ref{type}, $\mathcal{J}$ naturally induces a $G$
equivariant orthogonal map with square $-1$ on the $G$-invariant
vector bundle
$$(\g_M\oplus df)| _{f^{-1}(0)}^\perp  /
(\g_M \oplus df)|_{f^{-1}(0)}.$$ Let $\widetilde{\mathcal{J}} \colon
TM_0 \oplus T^* M_0 \to TM_0 \oplus T^*M_0$ be the induced
generalized almost complex structure on $M_0$.

Let $L \subset T_\C M \oplus T_\C ^*M$ and $\widetilde{L} \subset
T_\C M_0 \oplus T_\C ^* M_0$ be the $\sqrt{-1}$ eigenbundles of
$\mathcal{J}$ and  $\widetilde{\mathcal{J}}$, respectively. By the
definition of generalized complex structure, $L$ is closed under the
Courant bracket. By Lemma~\ref{le:closed1}, the image of $L \cap
df^\perp$ in $T_\C (f^{-1}(0)) \oplus T_\C ^*(f^{-1}(0))$ is also
closed under the Courant bracket. Since $L \cap df^\perp = L \cap
(\g_M \oplus df)^\perp$, by Lemma~\ref{le:closed2} $\widetilde{L}$,
which is its image  in $T_\C M_0 \oplus T_\C ^* M_0$, is also
closed.

The last statement is a direct consequence of Lemma \ref{type}.

This proves the case  $a = 0$.

%%S added $\mu_a$ and $M_a$.

For $a \neq 0$, let $\omega$ be the Kirrilov-Kostant symplectic form
on the co-adjoint orbit $\cO_{-a}$ and let $\mathcal{J}_{\omega}$ be
the induced generalized complex structure. Then $(\cO_{-a},
\mathcal{J}_{\omega})$ is a generalized complex manifold of type
$0$, and inclusion is a generalized moment map for the co-adjoint
$G$ action. Hence, $(M \times \cO_{-a}, \J \times
\mathcal{J}_{\omega})$ is a generalized complex manifold, $\type(\J
\times \mathcal{J}_{\omega})_{(m,b)} = \type(\J)_m$, and $\mu_a(x,v)
= \mu(x) - v$ is a generalized moment map for the diagonal action of
$G$ on $M \times \cO_{-a}$. Since it is easy to see that $M_a$ can
be identified with $\mu_{a}^{-1}(0) / G$, the result follows  from
the case $a = 0$.
\end{proof}

We  now find the generalized complex quotients for our basic
examples from Example~\ref{ex1} and \ref{ex2}.

\begin{example} \label{ex3}\
\begin{itemize}
\item [a)]  The generalized complex quotient of
the Hamiltonian generalized complex manifold associated to a
Hamiltonian symplectic manifold is the generalized complex manifold
associated to the symplectic quotient.
\item [b)] Since there is no generalized moment map in the
complex case, there is no generalized complex quotient.
\item [c)] As in the symplectic case, the generalized complex
quotient for the diagonal action on the product of two generalized
complex manifolds is not the product of the quotients.
\item [d)]
Let a compact Lie group $G$ act on a generalized complex  manifold
$(M,\J)$ with a real generalized moment map $ f \colon M \to
\g^*_\C$. Let $\J'$ be the $B$-transform of $\J$, where $B \in
\Omega^2(M)$ is closed, invariant and satisfies $\iota_{\xi_M} B =
0$ for all $\xi \in \g$. Then $B$ descends to a closed two-form
$\widetilde{B}$ on $M_a = f^{-1}(\cO_a)/G$ for any regular $a \in
\g^*$.  The generalized complex structure $\widetilde{\J}'$ for  the
generalized complex quotient of $(M,\J',\mu )$ at $a$ is the
$\widetilde{B}$ transform of the generalized complex quotient
$\widetilde{\J}$ of $(M,\J,\mu)$ at $a$. Moreover, if $B = d\gamma$,
where $\gamma \in \Omega^1(M)$ is closed, invariant and satisfies
$\iota_{\xi_M} \gamma$ for all $\xi \in \g$, then then $\gamma$
descends to a one-form $\widetilde{\gamma}$ on $M_a$, so that $d
\widetilde{\gamma} = \widetilde{B}$. Hence, in this case,
$\widetilde{\J}'$ is an exact $B$-transform of $\widetilde{\J}.$
\end{itemize}
\end{example}

Our first main result is now very easy to prove.

\begin{proposition} \label{Complex Reduction}
Let a compact Lie group $G$ act on a generalized complex manifold
$(M,\mathcal{J})$ with generalized moment map $\mu = f + \sqrt{-1} g
\colon M \to \g_\C^*$, where $f$ and $g$ are real. Let $\cO_a$ be
the co-adjoint orbit through $a \in \g^*$. If $G$ acts freely on
$f^{-1}(\cO_a)$, the generalized complex quotient $M_a$ inherits a
generalized complex structure $\widetilde{\J}$, which is natural up
to an exact $B$-transform. Moreover, for all $m \in f^{-1}(\cO_a)$,
$$\type(\widetilde{\mathcal{J}})_{[m]} = \type(\mathcal{J})_m.$$
\end{proposition}

\begin{proof}
By restricting to a neighborhood of $f^{-1}(\O_a)$, we may assume
that $G$ acts freely. Choose a connection $\theta \in
\Omega^1(M,\g)$, and let $B = d(h,\theta)$. Then by
Lemma~\ref{B-transform}, the $B$-transform $\J'$ of $\J$ is an
invariant generalized complex structure with real  moment map $f$.
Hence, by Lemma~\ref{complex reduction}, $\J'$ descends to a natural
generalized complex structure $\widetilde{\J}$ on $M_a$.

If we choose a different connection  $\widehat{\theta} \in
\Omega^1(M,\g)$ the resulting generalized complex structure
$\widehat{\J}'$ is the $d\gamma$ transform of $\J'$, where $\gamma =
( \theta' - \theta,h)$. Since $\gamma$ is invariant and
$\iota_{\xi_M} \gamma = 0$ for all $\xi \in \g$, the resulting
generalized complex quotient will be  an exact $B$-transform of
$\widetilde{\J}$. (See the last part  of Example~\ref{ex3}.)
\end{proof}

\section{Generalized  K\"{a}hler structure}\label{generalized kahler
structure}

%Added "n ordered"
A {\bf generalized K\"{a}hler structure} on a vector space $V$
consists of an ordered pair $(\J_1,\J_2)$ of commuting generalized
complex structures on $V$ so that
%%S small change
%$G :=-{\J}_1 {\J}_2$
$G =-{\J}_1 {\J}_2 \colon V \oplus V^* \to V \oplus V^*$ is a {\bf
positive definite metric}, by which we mean that
%%S MATH
$G^2 =1$,  $G$ is orthogonal, and $\langle G(w),w  \rangle > 0$ for
all non-zero $w \in V \oplus V^*$.
%%S MATH
Note that the first two conditions are automatically satisfied.

We will need the following lemma:

\begin{lemma}
%%S small change
Let $V$ be a vector space and let $W = V \oplus V^*$. Let $G \colon
W \to W$
%Let $V$ be a vector space and let $W = V \oplus V^*$. Let $G : W \to W $
 be a positive definite metric.
%%S slightly rewritten
Given an isotropic subspace $P$, define
$$
\widehat{W} = P^\perp \cap G(P)^\perp \subset W \qquad \mbox{and}
\qquad \widetilde{W} = P^\perp/P.
$$
The natural projection induces an isomorphism
\begin{equation}\label{iso}
\widehat{W} \iso \widetilde{W}.
\end{equation}
%Given an isotropic subspace $P$, the natural projection from
%$$\widehat{W} := P^\perp \cap G(P)^\perp \subset W
%$$
%to
%$$
%\widetilde{W} := P^\perp/P
%$$
%is an isomorphism.
\end{lemma}

\begin{proof}
For all non-zero $w \in P$, $\langle w, G(w) \rangle > 0$.
Therefore,  $G(P)^\perp \cap P = P^\perp \cap G(P) = \{0\}$. Since
$P$ is isotropic, this implies that $P \cap G(P) = \{0\}$, and hence
$\dim(\widehat{W}) =  \dim W  - 2 \dim P$. It also implies that the
projection above is an injection; the result now follows by a
dimension count.
\end{proof}

Let $(\J_1,\J_2)$ be a generalized K\"{a}hler structure on a vector
%%S small change + MATH change
space $V$, and let $G = - \J_1 \J_2$. If $P \subseteq W = V \oplus
V^*$
%space $V$, let $G = \J_1 \J_2$. If $P \subseteq W := V \oplus V^*$
%%S
%is any $\J_1$ invariant subspace, then by definition
is any $\J_1$ invariant subspace, then by since $\J_1 G = G \J_1$,
$P \oplus G(P)$ is $\J_1$ and $\J_2$ invariant and so we can define
the restrictions $\widehat{\J}_i$ of $\J_i$ to $\widehat{W} = W \cap
P^\perp \cap G(P)^\perp$. Clearly, $\widehat{\J}_1$ and
$\widehat{\J}_2$ are orthogonal and commute, $\widehat{\J}_1^2 =
%%S small change
%\widehat{\J}_2^2 = -1$, and $\widehat{G} :=  - \widehat{\J}_1
\widehat{\J}_2^2 = -1$, and $\widehat{G} =  - \widehat{\J}_1
\widehat{\J}_2$ is a positive definite metric on $\widehat{W}$. If
$P$ is also isotropic, then under the isomorphism (\ref{iso}) the
$\widehat{\J}_i$ induce maps $\widetilde{\J}_i \colon \widetilde{W}
\to \widetilde{W}$ satisfying the analogous conditions.
%%S moved and rewritten.
%It is easy to check that $\widetilde{\J}_1$ is exactly the map
%with the same name defined in the previous section.
%%S MATH
Moreover,  if $L_i$ is the $\sqrt{-1}$ eigenspace of $\J_i$, then
$L_i \cap \widehat{W}$ is the $\sqrt{-1}$ eigenspace of
$\widehat{\J}_i$, and the $\sqrt{-1}$ eigenspace of
$\widetilde{J}_i$ is its image under the isomorphism (\ref{iso}). As
in the previous section, if $P = (P \cap V) \oplus (P \cap V^*)$
then $(\widetilde{\J}_1,\widetilde{\J}_2)$ is a natural generalized
%%S small change
K\"{a}hler structure on $\widetilde{V} = (P \cap V^*)^0/(P \cap
%K\"{a}hler structure on $\widetilde{V} := (P \cap V^*)^0/(P \cap
V)$. It is easy to check that $\widetilde{\J}_1$ is the natural
complex structure on $\widetilde{V}$ defined in the beginning of the
previous section.

\begin{lemma}\label{type2}
Let $(\J_1,\J_2)$ be a generalized K\"{a}hler structure on a vector
space $V$. Consider $Q \subseteq V$ so that $\J_1(Q) \subseteq V^*$
and $P = Q \oplus \J_1(Q) \subset V \oplus V^*$ is isotropic.
%S
Let $(\widetilde{\J}_1,\widetilde{\J}_2)$ be the natural K\"ahler
strucutre on  $\widetilde{V}= \J_1(Q)^0/Q$.
%Define a generalized K\"{a}hler structure
%$(\widetilde{\J}_1,\widetilde{\J}_2)$ on $\widetilde{V}$ as in the
%paragraph above.
Then
$$\type(\widetilde{\J}_1) = \type(\J_1) \qquad \mbox{and}$$
$$\type(\widetilde{\J}_2) = \type(\J_2)-\dim(Q) +  2\dim(Q_\C\cap \pi(L_2)).$$
\end{lemma}

\begin{proof}
The first claim was proved in Lemma~\ref{type}.

We now turn to the second claim. Let $L_2$ and $\widetilde{L}_2$
denote the $\sqrt{-1}$ eigenspaces of $\J_2$ and $\widetilde{\J}_2$
respectively; let $\pi \colon V_\C \oplus V^*_\C \to V_\C$ and
$\widetilde{\pi} \colon \widetilde{V}_\C \oplus \widetilde{V}^*_\C
\to \widetilde{V}_\C$ be the natural projections.

Since $P$ is $\J_1$ invariant, $G(P) = \J_2(P)$, so $\J_2(P) \cap P
= \{0\}$. Moreover, by assumption, $\pi(P_\C) = \pi(Q_\C)$ and
$\dim(P) = 2 \dim(Q)$. Therefore, by Lemma \ref{pretype},
$\dim(\pi(L_2 \cap P_\C^\perp \cap \J_2(P_\C)^\perp)) = \dim(\pi(L_2
+ P_\C)) - \dim(P) = \dim(\pi(L_2) + Q_\C) - 2 \dim(Q) =
\dim(\pi(L_2)) - \dim(Q) - \dim(Q_\C \cap\pi( L_2))$. Moreover,
$\widetilde{\pi}(\widetilde{L}_2)$ is the projection of $\pi(L_2
\cap P_\C^\perp \cap \J_2(P_\C)^\perp) \subset \J_1(Q_\C)^0$ to
$\widetilde{V}_\C = \J_1(Q_\C)^0/Q_\C$, which reduces the dimension
by a further $\dim(\pi(L_2) \cap Q_\C)$. Finally,
$\dim(\widetilde{V}) = \dim(V) - 2 \dim(Q)$.

\end{proof}

A {\bf generalized K\"{a}hler structure} on a manifold $M$ is a pair
of commuting generalized complex structures $\J_1$ and $\J_2$ on $M$
so that $G = - \J_1 \J_2$ is a positive definite metric on $TM
\oplus T^*M$.

Again, we turn to \cite{Gua}  for our basic examples.

\begin{example} \label{Kahlerexamples} \cite{Gua} \
\begin{itemize}
\item [a)]
%%S slightly rewritten
%Let $(g, J,\omega)$ be a K\"{a}hler structure (in the usual sense)
%on a manifold $M$, i.e., a Riemannian structure $g$, a complex
%structure $J$, and a symplectic structure $\omega$ such that $g = -
%\omega J$.
%Then by Example \ref{ex1}, we have two
%generalized complex structures $\mathcal{J}_{\omega}$ and
%$\mathcal{J}_J$ which are induced by the symplectic structure
%$\omega$ and the complex structure $J$ respectively.
Let $(\omega,J)$ be a {\bf genuine K\"{a}hler structure} on a
manifold $M$, that is, a symplectic structure $\omega$ and a complex
structure $J$ which are {\bf compatible}, which means that $g = -
\omega J$ is a Riemannian metric. By Example \ref{ex1},  $\omega$
and $J$ induce generalized complex structures $\mathcal{J}_{\omega}$
and $\mathcal{J}_J$, respectively. Moreover, it is easy to see that
$\mathcal{J}_J$ and $\mathcal{J}_{\omega}$ commute, and that
\begin{equation} \label{metric}
G=-\mathcal{J}_\omega \mathcal{J}_{J}=\left(
\begin{matrix} 0 & g^{-1} \\g & 0 \end{matrix} \right)\end{equation}
is a positive definite metric on $TM \oplus T^*M$. Hence
$(\mathcal{J}_\omega, \mathcal{J}_J)$ is a generalized K\"{a}hler
structure on $M$.

\item [b)]
Let $(M,\J_{M,1}, \J_{M,2})$ and $(N,\J_{N,1},\J_{N,2})$ be
generalized K\"{a}hler manifolds, and define
$$\J_1= \left(
\begin{matrix} \J_{M,1} & 0
\\0 & \J_{N,1} \end{matrix} \right) \qquad \mbox{and} \qquad \J_2= \left(
\begin{matrix} \J_{M,2} & 0
\\0 & \J_{N,2} \end{matrix} \right).$$
Then $(M \times N, \J_1, \J_2)$ is a generalized K\"ahler manifold.
\end{itemize}
\end{example}

\begin{definition}
Let the compact Lie group $G$ with Lie algebra $\g$ act on a
manifold $M$. A {\bf generalized moment map}  for an invariant
generalized K\"{a}hler structure $(\J_1,\J_2)$ is a generalized
moment map for the generalized complex structure $\J_1$. (See
Definition~\ref{defmm}.)
\end{definition}

As before, let a compact Le group $G$ act on a generalized K\"ahler
manifold with generalized moment map $\mu = f + \newi h$, where $f$
and $h$ are real. Let $\cO_a$ be the coadjoint orbit through $a \in
\g^*$. If $G$ acts freely on $f^{-1}(\cO_a)$ then the {\bf
generalized K\"ahler quotient}
$$M_a =  f^{-1}(\cO_a)/ G$$ is a manifold.

\begin{example}\label{example Kahler 2} \
\begin{itemize}
\item [a)] If a compact Lie group $G$
acts on a K\"{a}hler manifold $(M,J,\omega)$ with moment map $\Phi$,
then $\Phi$ is  the generalized moment map for the $G$ action on
$(M,\mathcal{J}_J,\mathcal{J}_\omega)$.
\item [b)] If a compact Lie group $G$ acts on
two generalized K\"{a}hler manifold $(M,\J^M_1,\J^M_2)$ and
$(N,\J^N_1,\J^N_2)$ with moment maps $\mu^N$ and $\mu^N$, then
$\mu^M + \mu^N$ is a generalized moment map for the diagonal $G$
action on $(M \times N, \J^M_1 \times \J^N_1, \J^M_2\times\J^N_2)$.
\end{itemize}
\end{example}

We can now state our second main proposition:

\begin{proposition} \label{gkahler quotient}
Let a compact connected Lie group $G$ act on a generalized
K\"{a}hler manifold $(M,\J_1,\J_2)$ with generalized moment map $\mu
= f + \newi h \colon M \to \g_\C^*$. Let $\cO_a$ be the co-adjoint
orbit through $a \in \g^*$. If  $G$ acts freely on $f^{-1}(\cO_a)$,
then the sub-quotient $M_a = f^{-1}(\cO_a)/G$ naturally inherits a
generalized K\"{a}hler structure
$(\widetilde{\J}_1,\widetilde{\J}_2)$.

Moreover, let $\fk$ be the Lie algebra of the stabilizer $K$ of $a$,
%%S small change
and let $L_2$ be the $\sqrt{-1}$ eigenbundle of $\J_2$.
%%S killed
%and let $\pi \colon T_\C M \oplus T^*_\C M \rightarrow T_\C M$
%be the natural projection.
Then for all $m \in M$,
$$\type(\widetilde{\J}_1)_{[m]} = \type(\J_1)_m, \qquad \mbox{and}$$
$$\type(\widetilde{\J}_2)_{[m]} = \type(\J_2)_m - \frac{1}{2}\dim G
- \frac{1}{2}\dim(K) + 2 \dim (\fk_M \cap \pi(L_2))_m.$$
\end{proposition}

\begin{proof}
As before, we begin by assuming that $a = 0$.

By restricting to a neighborhood of $f^{-1}(0)$, we may assume that
%%S small change
%$G$ acts freely. By definition, $G := - \J_1 \J_2$ defines a
$G$ acts freely. By definition, $G = - \J_1 \J_2$ defines a
connection $\theta \in \Omega^1(M,\g^*)$  given by $\theta^\xi =
\frac{ G(\cdot,\xi_M)}{G(\xi_M,\xi_M)}$ for all $\xi \in \g$. By
Lemma~\ref{B-transform}, after apply an exact $B =  d(h,\theta)$
transform, we may assume that $h = 0$.

%%S small change
%As in the proof of Lemma~\ref{complex reduction}, by the
As in the proof of Proposition \ref{complex reduction}, by the
definition of generalized moment map $\g_M \oplus df$ is a $\J_1$
invariant subbundle, and $(\g_M \oplus df)|_{f^{-1}(0)}$ is
%%S small change
isotropic. Define $\widehat{W} =  (\g_M \oplus   df)^\perp \cap
%isotropic. Define $\widehat{W} :=  (\g_M \oplus   df)^\perp \cap
G(\g_M \oplus df)^\perp \subset TM \oplus T^*M$, and let
$\widehat{\J}_1$ and $\widehat{\J}_2$ be the restriction of $\J_1$
and $\J_2$ to $\widehat{W}$. Let
$(\widetilde{\J}_1,\widetilde{\J}_2)$ be generalized almost
K\"{a}hler structure on $M_0$ induced by $\widehat{\J}_1$ and
$\widehat{\J}_2$ under the restriction to $f^{-1}(0)$, isomorphism
from $\widehat{W}|_{f^{-1}(0)}$ to
%%S small change
%$(\g_M \oplus df)|_{f^{-1}(0)}^\perp \diagup (\g_M \oplus df)
$(\g_M \oplus df)|_{f^{-1}(0)}^\perp / (\g_M \oplus df)
|_{f^{-1}(0)},$ and the quotient map from $f^{-1}(0)$ to $M_0$.

%%S small change
%In Proposition \ref{complex reduction}, we checked that
In Lemma~\ref{complex reduction}, we checked that $\widetilde{\J}_1$
is a  generalized complex structure.

Let $L_2 \subset TM_\C \oplus T_\C^*M$, $\widehat{L}_2 \subset
%%S small change
%\widehat{W}_\C:=\widehat{W}\otimes \C$, and $\widetilde{L}_2 \subset
\widehat{W}_\C =\widehat{W}\otimes \C$, and $\widetilde{L}_2 \subset
T_\C M_0 \oplus T^*_\C M_0$ be $\sqrt{-1}$ eigenbundles of $\J_2$,
$\widehat{\J}_2$, and $\widetilde{\J}_2$, respectively. Since $\J_2$
is a generalized complex structure, $L_2$ is closed under the
Courant bracket. Since $G = -\J_1 \J_2$,  and $\J_1(\g_M) = df$,
$G(\g_M)  = \J_2(df)$ and $G(df) = \J_2(\g_M)$. Therefore
$\widehat{L}_2 = L_2 \cap \widehat{W}_\C = L_2 \cap \g_M^\perp \cap
df^\perp$. Therefore, by
Lemmas~\ref{le:closed1}~and~\ref{le:closed2}, $\widehat{L}_2$ is
also closed under the  Courant bracket. Moreover, since
$\widetilde{L}_2$ is the image of $\widehat{L}_2$ under the natural
restriction and quotient maps, it is also closed under Courant
bracket by the same lemmas. Therefore,
$(\widetilde{\J}_1,\widetilde{\J}_2)$ is a generalized K\"{a}hler
structure.

The formulas on types follow directly from Lemma \ref{type2}.

This proves the case  $a = 0$. For $ a \neq 0$, as in
Lemma~\ref{complex reduction} let $\omega$ be the Kirrilov-Kostant
symplectic form on $\cO{-a}$, and let $J$ be the natural invariant
complex structure which is compatible with $\omega$. Let
$(\mathcal{J}_{\omega},\mathcal{J}_J)$ be the induced generalized
K\"ahler structure on $\cO_{a}$; it has type $(0,\frac{1}{2} \dim
(\cO_{-a}))$ and inclusion is a generalized moment map for the
co-adjoint $G$ action on $\cO_{-a}$. Hence $(M \times O_{-a}, \J_1
\times \mathcal{J}_\omega, \J_2 \times \mathcal{J}_J)$ is a
generalized K\"ahler manifold,  $\type(\J_1 \times
\mathcal{J}_\omega)_{(m,b)} = \type(\J_1)_m$ and $\type(\J_2 \times
\mathcal{J}_J)_{(m,b)} = \type(\J_2)_m  + \frac{1}{2} \dim(O_a)$ for
all $m \in M$, and $\mu_a(x,v) = \mu(x)-v$ is a generalized moment
map for the diagonal action of $G$ on $M \times \cO_{-a}$. Finally,
it is easy to check that the intersection of $\g_{M \times O_{-a}}$
with the projection of the $\sqrt{-1}$ eigenbundle of $\J_2 \times
\mathcal{J}_J$ to $T_\C(M \times O_{-a})$ is isomorphic to $\fk_M
\cap \pi(L_2)$. Since $M_a$ can be identified with $\mu_{a}^{-1}(0)/
G$, the result follows from  case that $a = 0$.
\end{proof}

%% Totally rewritten
\begin{example} \label{Kahler quotient}
The generalized K\"ahler quotient of the Hamiltonian generalized
K\"ahler manifold associated to a Hamiltonian K\"ahler manifold is
the generalized K\"ahler manifold associated to the K\"ahler
quotient.
\end{example}

\begin{example} \label{hyper-kahler}Let $(M,g,I,J,K)$ be a hyper K\"{a}hler
structure, and let $\omega_I, \omega_J$ and $\omega_K$ be the
K\"{a}hler two forms that correspond to the complex structure $I,J$
and $K$ respectively. As shown in \cite{Gua}, we can construct a
generalized K\"{a}hler structure $(\mathcal{J}_1,\mathcal{J}_2)$ as
follows:
\begin{equation} \mathcal{J}_1 = \left( \begin{matrix} 1 & o
\\\omega_K & 1 \end{matrix} \right) \left(\begin{matrix} 0 &
-\dfrac{1}{2}(\omega_I^{-1}-\omega_J^{-1})\\\omega_I-\omega_J &0
\end{matrix} \right) \left( \begin{matrix} 1 & o
\\-\omega_K & 1 \end{matrix} \right) \end{equation}
\begin{equation} \mathcal{J}_2= \left( \begin{matrix} 1 &0
\\-\omega_K & 1 \end{matrix} \right) \left(\begin{matrix} 0 &
-\dfrac{1}{2}(\omega_I^{-1}+\omega_J^{-1})\\\omega_I+\omega_J &0
\end{matrix} \right) \left( \begin{matrix} 1 & 0
\\ \omega_K & 1 \end{matrix} \right) \end{equation}

Suppose there is a $S^1$ action on $M$ with the fundamental vector
field $X$ such that $\iota_X \omega_I =d\mu_I, \iota_X
\omega_J=d\mu_J,$ and $ \iota_X \omega_K=d\mu_K$ for some smooth
functions $\mu_I, \mu_K,$ and $\mu_J$. Set $f=\mu_I-\mu_J$. Then

 \begin{equation} \begin{split} \mathcal{J}_1 df &=
\left( \begin{matrix} 1 & 0
\\ \omega_K & 1 \end{matrix} \right) \left(\begin{matrix} 0 &
-\dfrac{1}{2}(\omega_I^{-1}-\omega_J^{-1})\\\omega_I-\omega_J &0
\end{matrix} \right) \left( \begin{matrix} 1 & 0
\\-\omega_K & 1 \end{matrix} \right) \left(\begin{matrix} 0\\df
\end{matrix}\right)\\ &=\left( \begin{matrix} 1 & 0
\\\omega_K & 1 \end{matrix} \right) \left(\begin{matrix} 0 &
-\dfrac{1}{2}(\omega_I^{-1}-\omega_J^{-1})\\\omega_I-\omega_J &0
\end{matrix} \right)\left(\begin{matrix} 0\\df
\end{matrix}\right)\\ &= \left( \begin{matrix} 1 & 0
\\\omega_K & 1 \end{matrix} \right)  \left(\begin{matrix} -X \\0
\end{matrix} \right) \\&=\left(\begin{matrix} -X  \\-d\mu_K
\end{matrix} \right) \end{split} \end{equation}

Thus  $f+\newi \mu_K$ is a generalized moment map for the circle
action on the generalized K\"{a}hler manifold $(M, J_1,J_2)$. Let
$\theta= \frac{g(\cdot , X)}{g(X,X)}$,  let $B=- d(\mu_K\theta)$,
and let $\mathcal{J}'_i$ be the $B$-transform of $\mathcal{J}_i$.
Then $(\mathcal{J}'_1, \mathcal{J}_2)$ is a generalized K\"{a}hler
structure which satisfies  $\mathcal{J}'_1 df=X$ on the level set
$f^{-1}(0)$. Assume that $S^1$ acts freely on $f^{-1}(0)$.
Proposition \ref{gkahler quotient} then asserts that there is a
reduced generalized K\"{a}hler structure
$(\widetilde{J}_1,\widetilde{J}_2)$ on the quotient
%%S small changes
$M_0 =f^{-1}(0)/ S^1$.
%$M_0:=f^{-1}(0)\diagup S^1$.

Assume in addition  that $0$ is a regular value for the map
$$\mu=(\mu_I,\mu_J,\mu_K): M \to R\oplus R \oplus R.$$
Then the hyper-K\"{a}hler quotient $S:= \left(\mu_I^{-1}(0)\cap
\mu_J^{-1}(0) \cap \mu_K^{-1}(0)\right)
%%S small change
%\diagup S^1$ is a codimension
/ S^1$ is a codimension two sub-manifold sitting inside $M_0$. We
have a natural inclusion map $i: S \to M_0$.

Let $(\widetilde{J}^S_1, \widetilde{J}^S_2)$ be the generalized
K\"{a}hler structure on $S$ induced by the quotient hyper-K\"{a}hler
structure on $S$. Let $\widetilde{L}^S_i$ be the $\sqrt{-1}$
eigenbundle of $\widetilde{J}^S_i$, and let $\widetilde{L}_i$ be the
$\sqrt{-1}$ eigenbundle of $\widetilde{J}_i$, $i=1,2$. Then we have
that $\widetilde{J}^S_i=\{X+ i^* \alpha  \mid X+\alpha \in
\widetilde{L}_i \cap \left(T_{\C}S \oplus T^*_{\C}M_0\right) \}$,
$i=1,2$. The submanifolds of a generalized complex manifold is
studied extensively in \cite{BB03}. Using their terminology,  we see
that $(S, \widetilde{J}^S_i)$ is exactly a generalized complex
submanifold of $(M_0, \widetilde{J}_i)$.

\end{example}

\section{Constructing bi-Hermitian structures}

In this section we are going to present a simple explicit
constructions  of bi-Hermitian structure on $\CP^N$, Hirzebruch
surfaces, $\C\P^2$ blown up at an arbitrary number of points, and
complex Grassmannians.

We will do this by constructing non-standard generalized K\"ahler
structures on these spaces. Since each of these manifolds can be
expressed as a symplectic quotient of $\C^n$, we start with the
standard K\"{a}hler structure on $\C^n$. Using the deformation
theory for generalized complex structures developed in \cite{Gua},
we deform this to another invariant generalized K\"{a}hler
structure. These techniques are particularly easy and explicit in
this very simple example; we do not need to resort to any global
analysis. Then we use the quotient construction we developed in
Section \ref{generalized kahler structure} to construct a
generalized K\"{a}hler structure on the quotient space which is not
the $B$-transform of a genuine K\"{a}hler structure,  although in
each case the first generalized complex structure is the one induced
from the standard symplectic structure. By the connection between
generalized K\"{a}hler structures and bi-Hermitian structures which
was established by Gualtieri \cite{Gua}, and which we explain below,
this induces a bi-Hermitian structure on each manifold.

\subsection{Review}

We begin with a brief review; all the material in this subsection,
with the exception of material specifically attributed to other
authors, was taken from \cite{Gua}.

\begin{definition} (\cite{AGG99})       A bi-Hermitian structure on a manifold
$M$ is a triple $(g,J_+,J_-)$, where $g$ is a Riemannian metric and
$J_{+}$ and $J_{-}$ are complex structures which are orthogonal
(with respect to $g$), induce the same orientation, and satisfy
$J_1(x) \neq \pm J_2(x)$ for some $x \in M$.
\end{definition}

Given a generalized K\"{a}hler manifold $(M,\J_1,\J_2)$, let $G = -
\J_1 \J_2$ be the associated positive definite  metric. Recall that
$G^2 = 1$, and let $C_+$ denote the $+1$ eigenspace of $G$. Since
$C_+ \subset TM \oplus T^*M $ is positive definite and $T^* M$ is
isotropic, the natural projection $\pi \colon C_+ \to TM$ is an
isomorphism. Therefore, $\langle \cdot, \cdot \rangle$ descends to a
Reimannian metric $g$ on $M$. Since $\J_1$ and $\J_2$ commute with
$G$, they both preserve $C_+$. Therefore, $\J_1$ and $\J_2$ descend
to almost complex structures $J_+$ and $J_-$ on $M$ which are
orthogonal (with respect to $g$). Combining Proposition 6.15 with
Remarks 6.13 and 6.14 in \cite{Gua} we see, (after possibly
replacing $J_-$ by $-J_-$)

\begin{proposition}\cite{Gua}\label{equivalence}
Given a  generalized K\"{a}hler manifold $(M,\J_1,\J_2)$, the above
construction defines a bi-Hermitian structure $(g,J_+,J_-)$ on $M$
exactly if
\begin{enumerate}
\item the generalized K\"{a}hler
structure is not the $B$-transform of a genuine K\"{a}hler structure
on $M$;
\item  At least one of the $\J_i$ has even type.
\end{enumerate}
\end{proposition}

Note that, in fact, if $\dim M = 4k$, then either both $\J_1$ and
$\J_2$ have odd type, or they both have even type. In contrast, if
$\dim M = 4k + 2$, then one must always have odd type whereas the
other has even type; therefore, the second condition is empty.

Let $\J$  be a generalized complex structure on a vector space $V$.
Let $L \subset V_\C \oplus V_\C ^*$ be the $\sqrt{-1}$ eigenspace of
$\J$. Since $L$ is maximal isotropic and $L \cap \overline{L} =
\{0\}$, we can (and will) use the  metric to identify $L^*$ with
$\overline{L}$.

Given $\epsilon \in \wedge^2 L^*$, define
$$L_\epsilon = \{Y + \iota_Y\epsilon \mid Y \in L \}.$$
Then $L_\epsilon$ is  maximal isotropic, and $L_\epsilon \cap
\overline{L_\epsilon} = \{0\}$ if and only if the endomorphism
\begin{equation} \label{real index zero condition}  A_{\epsilon}=\left( \begin{matrix} 1 &  \bar{\epsilon} \\
\epsilon & 1 \end{matrix} \right)  \colon L\oplus \overline{L}
\rightarrow L \oplus \overline{L} \end{equation} is invertible. If
it is invertible, there exists a unique generalized complex
structure $\J_\epsilon$ on $V$ whose $\sqrt{-1}$ eigenspace is
$L_\epsilon$. Note that $A_\epsilon$ is always invertible for
$\epsilon$ sufficiently small.

Now let $(\J_1, \J_2)$ be a generalized K\"{a}hler structure on $V$.
Let $L_1$ and $L_2$ denote the $\sqrt{-1}$ eigenspaces of $\J_1$ and
$\J_2$, respectively. Then $L_1 =\left(L_1 \cap L_2 \right) \oplus
\left( L_1 \cap \overline{L_2} \right)$ and $L_2 =\left(L_1 \cap L_2
\right) \oplus \left( \overline{L_1}\cap L_2  \right)$. Thus
$\epsilon \in C^\infty(\wedge^2 \overline{L_2})$ fixes $\J_1$ if and
only if $\epsilon$ takes $L_1 \cap L_2$ to $L_1 \cap
\overline{L_2}$,
%%S mathematical difference.  We  didn't  have condition below.
%%S Which is right?
i.e., if and only if $\epsilon$ is an element of
$C^\infty\left((\overline{L_1} \cap\overline{L_2}) \otimes (L_1 \cap
\overline{L_2})\right)$.

%Then $L_2$ decomposes as $L_2=L_2^+ \oplus L_2^-$, and
%$L_1$ decomposes as $L_1= L_2^+ \oplus \overline{L^-_2}$. Suppose
%$\epsilon \in C^{\infty}(\wedge^2 \overline{L_2})$. Then it is easy
%to see $\epsilon$ keeps $\mathcal{\J}_1$ fixed if and only if
%$\epsilon$ takes $L_2^+$ to $\overline{L^-_2}$, i.e., if and only if
%$\epsilon$ is an element of $\overline{L_2^+}\otimes
%\overline{L_2^-}$.

We are now ready to state the condition for $L_\epsilon$ to be
closed under the Courant bracket, as proved in \cite{LWP97},
following the presentation in \cite{Gua}.

We begin with two definitions. Although both can be defined more
generally for any Lie algebroid, we will only state them for the
case which interests us.

\begin{definition}
Let $L \subset T_\C M \oplus T^*_\C M$ be an isotropic subbundle
which is closed under the Courant bracket and let $\pi \colon L \to
T_\C M$ denote the natural projection. The {\bf Schouten bracket} is
the $\R$-bilinear map
$$ [\cdot,\cdot] \colon C^{\infty}(\wedge ^p L) \times C^{\infty}(\wedge ^q L)
\rightarrow C^{\infty}(\wedge ^{p+q-1} L)$$ which is characterized
by the following two formulas:
$$ [X_1 \wedge
\cdots \wedge X_p, Y_1 \wedge  \cdots \wedge
Y_q]=\sum_{i,j}(-1)^{i+j}[X_i, Y_j]\wedge X_1 \wedge \cdots\wedge
\widehat{X}_i \wedge\cdots \wedge\widehat{Y}_j \wedge\cdots \wedge
Y_q $$ for any $X_i$ and $Y_j$ in $C^{\infty}(L)$, and
$$[Y, f]=-[f,Y]=\pi(Y)f $$
for any $Y \in C^{\infty}(L)$ and $f \in C^{\infty}(M).$
\end{definition}
\begin{definition}
Let $L \subset T_\C M \oplus T^*_\C M$ be an isotropic subbundle
which is closed under the Courant bracket and let $\pi \colon L \to
T_\C M$ denote the natural projection. The {\bf Lie algebroid
derivative} is a first order linear differential operator from
$C^{\infty}(\wedge ^* L)$ to $C^{\infty}(\wedge ^{*+1} L)$ defined
by

\begin{equation}\begin{split}  d_L \sigma (X_0, \cdots, X_k) & =\sum _i (-1)^i
\pi(X_i) \sigma (X_0, \cdots, \widehat{X}_i, \cdots, X_k)\\& +\sum
_{i<j}(-1)^{i+j}\sigma ([X_i, X_j],X_0, \cdots, \widehat{X}_i,
\cdots, \widehat{X}_j,\cdots,X_k), \end{split}\end{equation}£¬ where
$\sigma \in C^{\infty}(\wedge ^k L^*)$ and $X_i \in C^{\infty}(L)$.
\end{definition}

\begin{example}
Lf $L$ is the $\sqrt{-1}$ eigenspace of the generalized complex
structure  $\J_J$ associated to a complex structure $J$, then $d_L$
is $\overline{\partial}$.
\end{example}

We will need the following special case of the theorem from
\cite{LWP97}.

 \begin{theorem} \label{closedness}
Let  $L \subset T_\C M \oplus T_\C^* M$ be a  maximal isotropic
subbundle so that $L \cap \overline{L} = \{0\}$ which is closed
under the Courant bracket. For any $\epsilon \in \wedge ^2
\overline{L}$,
%$$L_{\epsilon} :=\{Y+\iota_Y \epsilon : Y \in L \}$$
$$L_{\epsilon} =\{Y+\iota_Y \epsilon  \mid Y \in L \}$$
 is closed under Courant bracket if and only if $\epsilon$
 satisfies the Maurer-Cartan equation:
$$d_L \epsilon + \dfrac{1}{2} [\epsilon, \epsilon]=0 .$$
\end{theorem}

\subsection{Examples}

Now we are ready to turn to specific examples.

\begin{example}{Structures on $\C^n$}

We will begin by deforming the generalized K\"ahler structure $(
\mathcal{J}_{\omega},\mathcal{J}_J)$ on $\C^N$ which is induced by
the standard genuine K\"ahler structure $(\omega,J)$. (See Example
\ref{Kahlerexamples}). Note that while our ideas for deforming this
structure are taken entirely from \cite{Gua}, and we use many of his
observations, the deformation is much easier and more explicit in
this very simple case than in general.  In particular, while we use
many observations from that paper, our construction does not rely on
any of the deeper theorems.

Since
$$\overline{L_1} \cap \overline{L_2} =
\{Y-\newi \iota_Y \omega \mid Y\in T_{1,0}(M)\} \ \mbox{and} \ L_1
\cap \overline{L_2} =\{ Z+\newi \iota_Z \omega \mid Z\in T_{1,0}(M)
\},$$ for any global sections  $Y$ and $Z$ of $T_{1,0}(M)$
$$ \epsilon  = Y\wedge Z + \iota_Y \omega \wedge \iota_Z\omega $$
$$\qquad =   \frac{1}{2} (Y-\newi \iota_Y\omega ) \wedge (Z +\newi \iota_Z
\omega)- \frac{1}{2}(Z-\newi \iota_Z\omega ) \wedge (Y +\newi
\iota_Y \omega)$$ lies in $\C^\infty((\overline{L_1} \cap
\overline{L_2}) \otimes (\overline{L_1} \cap L_2)).$
%More generally,
%$$\epsilon = \beta +
%\omega^{-1} \beta \omega \in \C^\infty((\overline{L_1} \cap
%\overline{L_2}) \otimes (\overline{L_1} \cap L_2)) \qquad \forall\
%\beta \in \C^\infty(\wedge^2 T_{1,0}(M)).$$

If we restrict to  any open bounded subset $U$ of $\C^n$, then after
multiplying $\epsilon$ by  a sufficiently small positive number,
$A_\epsilon$ will be invertible. Thus $\epsilon$ deforms
$\mathcal{J}_J$ to a new generalized complex structure $\J_\epsilon$
on $U$ while keeping $\mathcal{J}_{\omega}$ fixed. Moreover,
$\pi(L_\epsilon)$ is spanned by $T_{1,0} \C^n, Y$, and $Z$. Thus,
$\type(\mathcal{J}_\epsilon)_z = n-2$ wherever $Y \wedge Z \neq 0$,
and $n$ at every other point.

Finally, the  following lemma gives a simple condition which
guarantees that $L_\epsilon$ is closed under the Courant bracket,
and hence that $(\mathcal{J}_{\omega},\J_\epsilon)$ is a generalized
K\"ahler structure.

\begin{lemma} \label{closedness2}
Assume that there exists a subset $I \subset (1,\ldots,n)$ so that
$$\epsilon  =
\sum_{i,j \in I}  F_{ij} (z) \frac{\partial}{\partial z_i} \wedge
\frac{\partial}{\partial z_j} + \sum_{i,j  \in I} F_{ij}(z)
d\overline{z_i} \wedge d\overline{z_j}.$$ If $F_{ij}$ is holomorphic
and $\dfrac{\partial{F_{ij}}}{\partial z_k} = 0$ for all $i, j$ and
$k \in I$, then $L_\epsilon$ is closed under the Courant bracket.
\end{lemma}

\begin{proof}
Since the Lie algebroid derivative $d_{L_J}$ is
$\overline{\partial}$ and  $F(z)$ is holomorphic,  $d_{L_J} \epsilon
= 0$. Hence, by Theorem~\ref{closedness}, the deformed generalized
almost K\"ahler structure will be a generalized K\"{a}hler structure
exactly if $[\epsilon,\epsilon] = 0$. This follows from the
calculation below:
$$
\left[F_{ij}(z) \frac{\partial}{\partial z_i} \wedge
\frac{\partial}{\partial z_j} , F_{kl}(z) \frac{\partial}{\partial
z_k} \wedge \frac{\partial}{\partial z_l} \right]
$$
$$
= F_{kl} \( -\frac{\partial F_{ij}}{\partial z_k}
\frac{\partial}{\partial z_l} + \frac{\partial F_{ij}}{\partial z_l}
\frac{\partial}{\partial z_k} \) \wedge \frac{\partial}{\partial
z_i} \wedge \frac{\partial}{\partial z_j} + F_{ij} \(-
\frac{\partial F_{kl}}{\partial z_i} \frac{\partial}{\partial z_j} +
\frac{\partial F_{kl}}{\partial z_j} \frac{\partial}{\partial z_i}
\) \wedge \frac{\partial}{\partial z_k} \wedge
\frac{\partial}{\partial z_l} = 0.
$$
Similarly,
$$
\left[F_{ij}(z) \frac{\partial}{\partial z_i} \wedge
\frac{\partial}{\partial z_j} , F_{kl} d\overline{z_k}\wedge
d\overline{z_l} \right]
$$
$$=
F_{ij} \( \frac{\partial F_{kl}}{\partial z_j}
\frac{\partial}{\partial z_i} - \frac{\partial F_{kl}}{\partial z_i}
\frac{\partial}{\partial z_j} \) \wedge d\overline{z_k}\wedge
d\overline{z_l} = 0 $$ Finally,
$$
\left[F_{ij} d\overline{z_i}\wedge d\overline{z_j}, F_{kl}
d\overline{z_k}\wedge d\overline{z_l} \right] = 0.
$$
\end{proof}

\end{example}

Suppose that a compact Lie group $G$ acts on $(\C^n,\omega,J)$ with
proper moment map $\Phi  \colon \C^n \to \g^*$. Consider $a \in
\g^*$ so that $G$ acts freely on $\Phi^{-1}(\cO_a)$; let $M_a =
\Phi^{-1}(\cO_a)$ denote the symplectic quotient.

Assume that there exists a subset $I \subset (1,\ldots,n)$ so that
$$\epsilon  =
\sum_{i,j \in I}  F_{ij} (z) \frac{\partial}{\partial z_i} \wedge
\frac{\partial}{\partial z_j} + \sum_{i,j  \in I} F_{ij}(z)
d\overline{z_i} \wedge d\overline{z_j}.$$ Moreover, assume that
$F_{ij}$ is holomorphic and $\dfrac{\partial{F_{ij}}}{\partial z_k}
= 0$ for all $i, j$ and $k$ in $I$. Since $\cO_a$ is bounded, by
multiplying $\epsilon$ by a sufficiently small constant we may
assume that $A_\epsilon$ is invertible on $\cO_a$. Then, applying
Lemma~\ref{closedness2}, $(\J_\omega,\J_\epsilon)$ is an invariant
generalized K\"ahler structure with generalized moment map $\Phi$.
Hence, by Proposition~\ref{gkahler quotient}, there is a natural
generalized K\"{a}hler Structure $(\widetilde{\J}_\omega,
\widetilde{\J}_\epsilon)$ on the  symplectic quotient $M_a$.
Moreover, $\widetilde{\J}_\omega$ has type $0$; in fact, it is the
generalized complex structure associated to the usual symplectic
structure on $M_a$. Hence, condition (ii) of
Proposition~\ref{equivalence} is automatically satisfied. So, by
Proposition~\ref{equivalence}, $(\widetilde{\J}_\omega,
\widetilde{\J}_\epsilon)$ will induce a bi-Hermitian structure on
the reduced space as long as it is not the $B$-transform of a
genuine K\"{a}hler structure.  To check this, it is enough to check
that $\type(\widetilde{\J}_\epsilon)_{[z]} \neq \frac{1}{2} \dim
M_a$ for at least some $[z] \in  M_a$. Since $\type(\J_\epsilon)_z <
N$ for generic $z \in \C^n$, by Proposition~\ref{gkahler quotient},
it is enough to check that $\fk_{\C^n} \cap \pi(L_\epsilon) = \{0\}$
at generic points.

\begin{example} ({\bf  $\CP^N$ for $N \geq 2$})

We now construct a bi-Hermitian structure on $\CP^N$ for $N \geq 2$.

Let $S^1$ act on $\C^{N+1}$ via $$\lambda \cdot (z_0,\ldots,z_N) =
(\lambda z_0,\ldots, \lambda z_N).$$ Note that this action preserves
the K\"{a}hler structure $(\omega,J)$. Moreover, $$\Phi(z) = \sum_i
\frac{1}{2}|z_i|^2$$ is a moment map, $S^1$ acts freely on
$\Phi^{-1}(1)$, and the reduced space $M_1 = \Phi^{-1}(1)/S^1$ is
$\CP^N$.

Let $$\epsilon =z_0^2 \frac{\partial}{\partial z_1} \wedge
\frac{\partial}{\partial z_2} + z_0^2 d\overline{z}_1 \wedge
d\overline{z}_2.$$ After multiplying $\epsilon$ by a sufficiently
small positive constant, $A_\epsilon$ is invertible, so $\epsilon$
deforms $(\mathcal{J}_\omega,\mathcal{J}_{J})$ to a new generalized
almost K\"{a}hler structure $(
\mathcal{J}_\omega,\mathcal{J}_\epsilon)$ on $\C^N$, so that
$\type(\J_\epsilon)_z = N+1$ if $z_0 = 0$ and is $N - 1$ otherwise.
Since  $z_0^2$ is holomorphic and  $\dfrac{\partial z_0^2}{\partial
z_1}=\dfrac{\partial z_0^2}{\partial z_2}=0$, by
Lemma~\ref{closedness2} is in fact
$(\mathcal{J}_\omega,\J_\epsilon)$ is a generalize K\"{a}hler
structure.

Since $\epsilon$ is $S^1$ invariant,
$(\mathcal{J}_\omega,\J_\epsilon)$ is also invariant. Hence, by
Proposition \ref{gkahler quotient}, there is a natural generalized
K\"{a}hler structure $(\widetilde{\mathcal{J}}_\omega,
\widetilde{\J}_\epsilon)$ on the quotient space $\CP^N =
\Phi^{-1}(1) / S^1$.

Note that  fundamental vector generated by the action is
$$X =\dfrac{\newi}{2}\sum_i \left( z_i \dfrac{\partial }{\partial z_i}
-\overline{z}_i \dfrac{\partial }{\partial \overline{z}_i} \right)
\,,$$ and that $X$ does not lie in $\pi(L_\epsilon)$ at any point of
$\C^{N+1}$, where $L_\epsilon$ is the $\sqrt{-1}$ eigenbundle of
$\mathcal{J}_\epsilon$. It follows immediately from Proposition
\ref{gkahler quotient} that
$\type(\widetilde{\mathcal{J}_\omega})_{[z]} = 0$ for all $[z] \in
\CP^N$, whereas $\type(\widetilde{\mathcal{J}}_\epsilon)_{[z]} = N$
if $z_0=0$, otherwise
$\type(\widetilde{\mathcal{J}}_\epsilon)_{[z]}=N-2$. By Proposition
\ref{equivalence} $(\widetilde{\mathcal{J}}_{\omega},
\widetilde{\mathcal{J}}_\epsilon)$ gives us a bi-Hermitian structure
on $\CP^N$.

In the case of $N=2$, the above construction actually gives us a
$SU(2)$-invariant bi-Hermitian structure.  Note that the standard
action of $SU(2)$ on $C^2$ can be extended to $C^{3}$ by letting
$SU(2)$ act on the first component trivially. This $SU(2)$ action
commutes with the standard $S^1$ action on $C^3$ and therefore
descends to a $SU(2)$ action on $\CP^2$. Since both $\omega$ and
$\epsilon$ are $SU(2)$-invariant, the deformed generalized
K\"{a}hler pair $(\mathcal{J}_{\omega}, \J_\epsilon)$ is
$SU(2)$-invariant as well. Since the $SU(2)$ action on $C^{3}$
commutes with the standard $S^1$ action, we conclude that the
quotient generalized K\"{a}hler structure must be $SU(2)$-invariant.
\end{example}

\begin{example}{\bf Toric varieties}

We will now construct bi-Hermitian structures on many, but not all,
toric varieties, including  all Hirzebruch surfaces and the blow up
of $\CP^2$ at arbitrarily many points.

Let an $n$ dimensional torus $T$ with Lie algebra $\t$ act on a
compact symplectic manifold $(M,\omega)$ with moment map $\Psi
\colon: M \to \t^*$. Let $\Delta \subset \t^*$ be the moment
polytope. Let $\eta_1, \ldots, \eta_N  \in \t$ be the primitive
outward normals to the facets of $\Delta$. Define $p \colon \R^N \to
\t$ by $p(e_i) = \eta_i$. Let $K$ be the kernel of the associated
map from $(S^1)^N $ to $\t$. Let $K$ act on $\C^N$ via  its
inclusion into $(S^1)^N$; let $\Phi \colon  \C^N \to \fk^*$ denote
the resulting moment map. There exists some $\xi \in \fk^*$ so that
$K$ acts freely on $\Phi^{-1}(\xi)$ and $M$ is equivariantly
symplectomorphic to the reduced space
$$M_\xi = \Phi^{-1}(\xi)/K.$$

Now assume that there exists $\alpha \in  \t^*$ so that
$\alpha(\eta_1) = \alpha(\eta_2) = -1 $, but $\alpha(\eta_i) \geq 0$
for all other $i$. Note that this condition is not satisfied for all
toric symplectic manifolds, even in two dimensions. For example, it
is not satisfied for $\CP^1 \times \CP^1$ blown up at the four fixed
points. On the other hand, it is satisfied in many cases, including
Hirzebruch surfaces and $\CP^2$ blown up at an arbitrary number of
points as long as one picks those points carefully, for example,
blow up in a sequence of points so that each point lies on
$[0,z_1,z_2]$.

%Consider the exact
%sequence:
%$$0 \mapsto \fk \mapsto \R^N \stackrel{p}{\mapsto} \t \mapsto 0$$
%and the dual sequence:
%$$0 \mapsto \t^* \stackrel{p^*}{\mapsto} (\R^N)^* \mapsto \fk^* \mapsto 0.$$

Since  $\alpha(\eta_i) \geq 0$ for all $i \geq 3$, we may define
$$
\epsilon = \( \prod_{i \geq 3} z_i^{\alpha(\eta_i)} \)
\frac{\partial}{\partial z_1} \wedge \frac{\partial}{\partial z_2} +
\( \prod_{i \geq 3} z_i^{\alpha(\eta_i)} \) d\overline{z_1}\wedge
d\overline{z_2}.$$ After multiplying $\epsilon$ by a sufficiently
small positve constant, $A_\epsilon$ is invertible, so $\epsilon$
deforms $(\J_\omega,\J_J)$ to a new generalized almost K\"{a}hler
structure $(\J_\omega,\J_\epsilon)$ on $\C^N$ so that
$\type(\J_\epsilon)_z = N -2 $ if $z_j \neq 0$ for all $j$ so that
$\alpha(\eta_j) > 0$, and otherwise is equal to $N$. Since
$=\prod_{i \geq 3} z_i^{\alpha(\eta_i)}$ is holomorphic and
annihilated by $\dfrac{\partial}{\partial z_1}$ and
$\dfrac{\partial}{\partial z_2}$, by Lemma~\ref{closedness2},
$(\J_\omega,\J_\epsilon)$ is a generalized K\"{a}hler structure. By
construction, $\epsilon$ is $K$-invariant, so
$(\J_\omega,\J_\epsilon)$ descends to a  generalized K\"{a}hler
structure $(\widetilde{\J}_\omega,\widetilde{\J}_\epsilon)$ on the
reduced space $M_\xi$. Finally, for any $\beta \in \fk$, let
$\beta_i$ denote the $i$'th coordinate of its natural inclusion into
$\R^N$. Then $$\beta_{\C^N} = \sum_i \dfrac{\newi \beta_i }{2}
\left( z_i \dfrac{\partial }{\partial z_i} -\overline{z}_i
\dfrac{\partial }{\partial \overline{z}_i} \right) \,.$$ Since our
assumptions rule out $\eta_1 = - \eta_2$, $\beta_i \neq 0$ for some
$i$ which is not  $1$ or $2$. Hence, $\beta_{\C^N} \not\in
\pi(L_\epsilon)$. Thus, we get a bi-Hermitian structure.
\end{example}

\begin{example}{\bf Grassmannians}

Consider the natural action of  $G =U(n)$ on $M = \C^n \otimes \C^m$
with the moment map $\Phi \colon M \to \g^*$ given by
\begin{equation} \Phi(z) = \sum_{j = 1}^m \left( \begin{matrix}
z_{j1} \overline{z_{j1}} & z_{j1} \overline{z_{j2}} & \cdots & z_{j1} \overline{z_{jn}} \\
z_{j2} \overline{z_{j1}} & z_{j2} \overline{z_{j2}} & \cdots & z_{j2} \overline{z_{jn}} \\
\cdots&\cdots&\cdots &\cdots&\\
z_{jn} \overline{z_{j1}} & z_{jn} \overline{z_{j2}} & \cdots & z_{jn} \overline{z_{jn}} \\
\end{matrix} \right). \end{equation}
Here, we have labeled the coordinates
$z_{11}\ldots,z_{1n},\ldots,z_{m1},\ldots,z_{mn}$, and identified
$\g^*$ with $n \times n$ matrices $A$ such that $A =
\overline{A}^t$. Let $\cO_I$ be the coadjoint orbit containing the
identity matrix $I$. Note that $G$ acts freely on
$\Phi^{-1}(\cO_I)$, and furthermore that the reduced space
$\Phi^{-1}(\cO_I)/G$ is the Grassmannian of $n$ planes in $\C^m$.

Of course, $\sum_{i=1}^n z_{i,1} \overline{z_{i,1}}$ is $U(n)$
invariant. For the same reason, $\sum_{i=1}^n z_{i,1} d
\overline{z_{i,2}}$  and $\sum_{i=1}^n z_{i,1} d \overline{z_{i,3}}$
are $U(n)$ invariant. Therefore,
$$\epsilon =
\(\sum_{i=1}^n z_{i,1} \frac{\partial}{\partial z_{i,2}} \) \wedge
\(\sum_{i=1}^n z_{i,1} \frac{\partial}{\partial z_{i,3}} \)  +
\(\sum_{i=1}^n z_{i,1} d \overline{z_{i,2}}\) \wedge \(\sum_{i=1}^n
z_{i,1} d \overline{z_{i,2}}\)$$ is also $U(n)$ invariant.

After multiplying $\epsilon$ by a sufficiently small positive
constant $A_\epsilon$ is invertible, so it deforms
$(\J_\omega,\J_J)$ to a new generalized almost K\"{a}hler structure
$(\J_\omega,\J_\epsilon)$ so that $\type(\J_\epsilon)_z = nm - 2$
unless $z_{i,1} = 0$ for all $i$, in which case it is $nm.$ By
Lemma~\ref{closedness2}, $(\J_\omega,\J_\epsilon)$ is in fact a
generalized K\"{a}hler structure.

Moreover, it is easy to see that $\g_M \cap L_\epsilon  = \{0\}$.
Therefore, this gives rise to a bi-Hermitian structure on the
Grassmannian.

\end{example}

\end{document}